\title{Infinite Paley graphs}

\author{
Gareth~A.~Jones\\
}

\documentclass[12pt]{article}
\usepackage{color}
\usepackage[latin1]{inputenc}
\usepackage{a4wide}
\usepackage{latexsym}
\usepackage{amsmath}
\usepackage{amssymb}
\usepackage{amsfonts}
\usepackage{tikz}
\newtheorem{thm}{Theorem}[section]

\newcommand{\F}{\mathbb F}

\date{}
\begin{document}

\maketitle

\begin{abstract}
\noindent Infinite analogues of the Paley graphs are constructed, based on uncountably many infinite but locally finite fields. Weil's estimate for character sums shows that they are all isomorphic to the random or universal graph of Erd\H os, R\'enyi and Rado. Automorphism groups and connections with model theory are considered.
\end{abstract}


\medskip


\section{Introduction} In 1963 Erd\H os and R\'enyi~\cite{ER} described two constructions of graphs which have subsequently become  well-known and well-understood parts of the landscape of graph theory. One construction gave a countably infinite family of finite graphs, defined deterministically, which later became known as the Paley graphs $P(q)$. The other gave a single countably infinite graph $R$ (or more precisely an uncountable family of mutually isomorphic countably infinite graphs), defined randomly and later variously named after Erd\H os, R\'enyi and Rado, who gave an alternative construction in~\cite{Rad} the following year. It is perhaps surprising that in the following half-century and more, a strong connection between these very different graphs seems to have received little notice, except in the world of model theory (see~\cite[Examples~1.3.6 and 1.8.3]{MS}), though there are hints to be found in papers such as~\cite{BEH, BT}. Perhaps this lacuna is less surprising when one realises that an essential ingredient in this connection comes from algebraic geometry, namely Weil's estimate for character-sums, used in his proof of the Riemann hypothesis for curves over finite fields.

The first aim of this paper is give a more combinatorial explanation of this connection by constructing, for each odd prime $p$, infinite analogues of the Paley graphs, defined over uncountably many locally finite fields of characteristic $p$, and its second aim is to show that these graphs are all isomorphic to $R$. The finite and infinite Paley graphs are described in \S\ref{Paley} and \S\ref{InfPaley}, and $R$ is described in \S\ref{Random}. The isomorphism is proved in \S\ref{Proof}, with remarks on the proof in \S\ref{Remarks}. The automorphism groups of these finite and infinite graphs are compared in \S\ref{Auto}, and the construction and the isomorphism with $R$ are extended in \S\ref{GenPaley} to the generalised Paley graphs introduced by Lim and Praeger in~\cite{LP}. The paper~\cite{ER} is revisited in \S\ref{Symm}.


\section{The Paley graphs and their inclusions}\label{Paley}

For each prime power $q=p^e\equiv 1$ mod~$(4)$ the {\sl Paley graph\/} $P(q)$ has as its vertex set the field ${\mathbb F}_q$ of $q$ elements, with vertices $x$ and $y$ adjacent if and only if $x-y$ is a quadratic residue (non-zero square) in ${\mathbb F}_q$. It is an  undirected strongly regular graph with parameters $v=q$ (the number of vertices), $k=(q-1)/2$ (their common valency), $\lambda=(q-5)/4$ and $\mu=(q-1)/4$ (the number of common neighbours of two adjacent or non-adjacent vertices).
See~\cite{Bro} for further basic properties of the Paley graphs.

These graphs were introduced, not (as is often asserted) by Paley~\cite{Pal} in 1933, but in the case $e=1$ in 1962 by Sachs~\cite{Sac}, as examples of self-complementary graphs, and in the general case $e\ge 1$ in 1963 by Erd\H os and R\'enyi~\cite[\S 1]{ER}, as part of their study of asymmetric graphs. Neither paper attached a name to these graphs; they appear to have been named around 1970, no doubt by analogy with Paley designs (see~\cite{Dem}) which, like the orthogonal matrices constructed by Paley in~\cite{Pal}, are based on properties of quadratic residues in finite fields. See~\cite{Jon} for a discussion of the history of these graphs.

A finite field ${\mathbb F}_{q}$ is a subfield of ${\mathbb F}_{q'}$ if and only if $q'$ is a power $q^e$ of $q$, in which case the subfield is unique, consisting of the solutions of the equation $x^q=x$. In this case, if $q\equiv 1$ mod~$(4)$ then $q'\equiv 1$ mod~$(4)$, so we have Paley graphs $P(q)$ and $P(q')$. Clearly each quadratic residue in ${\mathbb F}_q$ is also a quadratic residue in ${\mathbb F}_{q'}$, so $P(q)$ is a subgraph of $P(q')$. If $e$ is even then every element of ${\mathbb F}_q$ has a square root in ${\mathbb F}_{q'}$ (in fact, in the quadratic subfield ${\mathbb F}_{q^2}\subseteq{\mathbb F}_{q'}$) so the subgraph of $P(q')$ induced by $P(q)$ is a complete graph $K_q$. However, if $e$ is odd then an element of ${\mathbb F}_q$ has a square root in ${\mathbb F}_{q'}$ if and only if it has one in ${\mathbb F}_q$, so the induced subgraph is simply $P(q)$, that is, $P(q)$ is a full subgraph of $P(q')$.


\section{Infinite Paley graphs}\label{InfPaley}

Let $E$ be any set of odd integers which is closed under taking divisors and least common multiples. For any prime $p\equiv 1$ mod~$(4)$ let
\[{\mathbb F}_{p^E}:=\bigcup_{e\in E}{\mathbb F}_{p^e},\]
the direct limit of the direct system of fields ${\mathbb F}_{p^e}$ for $e\in E$ and inclusions between them. This is a field, infinite if and only if $E$ is, and locally finite in the sense that each finite subset is contained in a finite subfield. The finite subfields of ${\mathbb F}_{p^E}$ are just the fields ${\mathbb F}_{p^e}$ for $e\in E$, so distinct sets $E$ determine distinct (and non-isomorphic) fields ${\mathbb F}_{p^E}$. There are uncountably many sets $E$ satisfying the above conditions (consider, for example, the set of integers $e$ whose prime factors all belong to a given set of odd primes), so for each $p$ we obtain uncountably many non-isomorphic fields ${\mathbb F}_{p^E}$.

Now let us define
\[P(p^E)=\bigcup_{e\in E}P(p^e),\]
the direct limit of the Paley graphs $P(p^e), e\in E$, with respect to the embeddings described above. By our earlier remarks, each $P(p^e)$ is a full subgraph of $P(p^E)$. If $E$ is finite then $E$ is just the set of all divisors of $l:={\rm lcm}(E)$, so that $P(p^E)$ is just another Paley graph $P(p^l)$. We will therefore assume from now on that $E$ is infinite, in which case we will call $P(p^E)$ an {\sl infinite Paley graph} (see~\cite[Example~1.8.3]{MS} for a similar construction by Macpherson and Steinhorn, though the exponents $e=2^i$ used there should be replaced with odd integers). In the same way one can construct infinite Paley graphs $P(p^{2E})$ for primes $p\equiv -1$ mod~$(4)$ as unions of Paley subgraphs $P(p^{2e})$ where $e$ is odd.

Despite the fact that they are constructed from uncountably many mutually non-isomorphic fields, these infinite Paley graphs $P(p^E)$ and $P(p^{2E})$ are all isomorphic to each other. In fact, we shall prove:

\begin{thm}\label{PaleyRthm}
Each infinite Paley graph $P(p^{rE})$ for $r=1,2$ is isomorphic to the random graph $R$.
\end{thm}


\section{The countable random graph}\label{Random}

The {\sl countable random graph}, or {\sl universal graph\/} $R$ was introduced by Erd\H os and R\'enyi~\cite[\S 3]{ER} in 1963 and Rado~\cite{Rad} in 1964. For details of its properties, see~\cite{Cam90, Cam97, Cam01} or~\cite[Section~9.6]{DM}. Theorem~\ref{PaleyRthm} should not be as surprising as it might at first appear, since in a sense we shall now explain `almost all' countably infinite graphs are isomorphic to $R$. 

Erd\H os and R\'enyi showed that if a graph $\Gamma$ has a countably infinite vertex set, and its edges are chosen randomly, then with probability $1$ it has the following property $U$: given any two disjoint finite sets $A$ and $B$ of vertices of $\Gamma$, there is a vertex which is a neighbour of each vertex in $A$ and a non-neighbour of each vertex in $B$. They used this to show that $\Gamma$ has a non-identity automorphism with probability $1$ (by contrast with the finite case, where a random graph of order $n$ has trivial automorphism group with probability approaching $0$ as $n\to\infty$). In fact, a similar argument shows that any two countably infinite graphs with property $U$ are isomorphic: one can construct an isomorphism between them by using $U$ to extend, by a back-and-forth argument, one vertex at a time, any isomorphism between finite induced subgraphs, such as a single vertex in each of them. (See, for example, \cite[Theorem~2.4.2]{Mar}, which in the language of model theory shows that the theory of graphs with property $U$ is satisfiable and $\aleph_0$-categorical, and hence complete and decidable.)  Thus any two graphs $\Gamma$ constructed randomly as above are isomorphic with probability $1$.

As a model of $R$ one can therefore take any countably infinite graph with property $U$. For instance, Rado~\cite{Rad} constructed a `universal graph', in which every countable graph is embedded as an induced subgraph, by using the vertex set $V=\mathbb N$ (including $0$), with vertices $x<y$ adjacent if and only if $2^x$ appears in the binary representation of $y$ as a sum of distinct powers of $2$; this easily implies property $U$.

The following well-known alternative model of $R$ imitates the construction of the (finite) Paley graphs. Let the vertex set $V$ be the (countably infinite) set of all primes $p\equiv 1$ mod~$(4)$, and define distinct vertices $p$ and $q$ to be adjacent if and only if $q$ is a quadratic residue mod~$(p)$, that is, the Legendre symbol $(\frac{q}{p})=1$. By quadratic reciprocity, which states that $(\frac{p}{q})(\frac{q}{p})=1$ for primes $p, q\equiv 1$ mod~$(4$), this is a symmetric relation, so it defines an undirected graph. To show that this graph has property $U$, given disjoint finite subsets $A$ and $B$ of $V$, for each prime $a\in A$ choose an integer $n_a$ which is a quadratic residue mod~$(a)$, and for each prime $b\in B$ choose an integer $n_b$ which is a non-residue mod~$(b)$. By the Chinese Remainder Theorem, the simultaneous congruences $n\equiv 1$ mod~$(4)$ and $n\equiv n_c$ mod~$(c)$ for all $c\in C:=A\cup B$ have a unique solution $n$ mod~$(d)$ where $d=4\prod_{c\in C}c$, and by a theorem of Dirichlet this congruence class contains a prime (infinitely many, in fact). This gives a vertex in $V$ adjacent to all the vertices $a\in A$ and to none of the vertices $b\in B$, as required.

\section{Proof of Theorem~\ref{PaleyRthm}}\label{Proof}

In order to prove Theorem~\ref{PaleyRthm} it is sufficient to prove that the infinite Paley graphs all have property $U$. To do this we will show that, given any two disjoint finite sets $A$ and $B$ of elements of a finite field ${\mathbb F}_q$ ($q$ odd), for all sufficiently large $e$ there is an element $x\in {\mathbb F}_{q^e}$ such that $x-a$ is a quadratic residue in ${\mathbb F}_{q^e}$ for all $a\in A$ and $x-b$ is a non-residue in ${\mathbb F}_{q^e}$ for all $b\in B$.

We will adapt an argument used by Blass, Exoo and Harary~\cite{BEH} to obtain a similar result concerning the family of Paley graphs $P(p)$ for primes $p\equiv 1$ mod~$(4)$. Given such subsets $A$ and $B$ of ${\mathbb F}_q$, let $S$ be the set of all $x\in{\mathbb F}_{q^e}$ satisfying the above condition. Let $C:=A\cup B$, let $n=|C|$, and let $\chi:{\mathbb F}_{q^e}\to{\mathbb C}$ be the quadratic residue character of ${\mathbb F}_{q^e}$, defined by $\chi(x)=1$, $-1$ or $0$ as $x$ is a quadratic residue, a non-residue or $0$. Note that $\chi(xy)=\chi(x)\chi(y)$ for all $x,y\in{\mathbb F}_{q^e}$.

For each $x\in{\mathbb F}_{q^e}\setminus C$ we have
\begin{equation}\label{abeqn}
\prod_{a\in A}(1+\chi(x-a)).\prod_{b\in B}(1-\chi(x-b))=
\begin{cases}
2^n & {\rm if} \; x\in S,\\
0 & {\rm otherwise}.
\end{cases}
\end{equation}
It follows that $S$ is non-empty if and only if
\[s:=\sum_{x\not\in C}\left(\prod_{a\in A}(1+\chi(x-a)).\prod_{b\in B}(1-\chi(x-b))\right)>0.\]
Summing over {\sl all\/} of ${\mathbb F}_{q^e}$ instead, let us define
\[t:=\sum_{x\in{\mathbb F}_{q^e}}\left(\prod_{a\in A}(1+\chi(x-a)).\prod_{b\in B}(1-\chi(x-b))\right).\]
Expanding the product on the right-hand side, we have
\[t=\sum_x 1+\sum_x\sum_a \chi(x-a)-\sum_x\sum_b \chi(x-b)+\cdots,\]
where the first term is $q^e$ and the second and third are $0$. To aid our consideration of the remaining terms, let us write $C=\{c_1,\ldots, c_n\}$. Then it follows from the above that
\[|t-q^e|\le|\sum_x\sum_{i_1<i_2}\chi(x-c_{i_1})\chi(x-c_{i_2})|+\cdots+
|\sum_x\sum_{i_1<\cdots<i_k}\chi(x-c_{i_1})\cdots\chi(x-c_{i_k})|+\cdots.\]
where $i_1,\ldots, i_k\in\{1,\ldots, n\}$ in each case. Weil's estimate~\cite{Wei} for character sums implies that
\begin{equation}\label{charsum}
\sum_x\chi(x-c_{i_1})\cdots\chi(x-c_{i_k})=O(q^{e/2})\quad\hbox{as}\;\;e\to\infty
\end{equation}
for each such $k$-tuple $(i_1,\ldots, i_k)$ (see Remark~1 for further explanation),
so it follows immediately  that
\[|t-q^e|=O(q^{e/2})\quad\hbox{as}\;\;e\to\infty.\]
Now
\[t-s=\sum_{x\in C}\left(\prod_{a\in A}(1+\chi(x-a)).\prod_{b\in B}(1-\chi(x-b))\right)\]
depends only on the sets $A$ and $B$, and not on $e$, so we have $s>0$ for all sufficiently large $e$, as required.


\section{Remarks on the proof}\label{Remarks}

\noindent{\bf 1.} Weil proved in~\cite{Wei} that if $\chi$ is a multiplicative character of order $d$ of a finite field ${\mathbb F}_q$ (one whose values are the $d$th roots of $1$ in $\mathbb C$), and $f(x)$ is a polynomial of degree $k$ over ${\mathbb F}_q$ not of the form $cg(x)^d$ for any $c\in{\mathbb F}_q$ and $g(x)\in{\mathbb F}_q[x]$, then
\begin{equation}\label{Weil}
\left|\sum_{x\in{\mathbb F}_q}\chi(f(x))\right|\le(k-1)\sqrt q.
\end{equation}
(See~\cite[p.~53]{Sch}, for example.) Replacing $q$ with $q^e$, taking $\chi$ to be the quadratic residue character, which has degree $d=2$, and taking $f(x)=(x-c_{i_1})\cdots(x-c_{i_k})$ we obtain the estimate (\ref{charsum}) used above. 

\medskip

\noindent{\bf 2.} The argument used to prove Theorem~\ref{PaleyRthm} in fact shows that
\[|S|=\frac{s}{2^n}\sim \frac{q^e}{2^n}\quad\hbox{as}\;\;e\to\infty,\; n\;\hbox{fixed},\]
which is what one would expect for Paley graphs on heuristic grounds, regarding adjacency or non-adjacency of vertices as independent events with equiprobable outcomes. Bollob\'as and Thomason~\cite{BT} have given a more precise estimate, equivalent in our notation to
\[\bigl||S|-\frac{q^e}{2^n}\bigr|\le\frac{1}{2}(n-2+2^{1-n})q^{e/2} +\frac{n}{2}.\]

\medskip

\noindent{\bf 3.} In~\cite{BEH} Blass, Exoo and Harary, working with the Paley graphs $P(p)$ for primes $p\equiv 1$ mod~$(4)$, needed to show that given any integer $n\ge 1$, if $p$ is sufficiently large then for any disjoint $n$-element sets $A$ and $B$ of vertices of $P(p)$ there is a vertex $x$ adjacent to every $a\in A$ and to no $b\in B$. Their argument (based on one for tournaments by Graham and Spencer~\cite{GS}) was similar to that used in \S\ref{Proof}, except that in place of Weil's character sum estimate for fields ${\mathbb F}_q$ they used one by Burgess~\cite{Bur}, that if $p$ is prime and $c_1,\ldots, c_k$ are distinct elements of ${\mathbb F}_p$, then
\[\bigl|\sum_{x\in{\mathbb F}_p}\chi(x-c_1)\ldots\chi(x-c_k)\bigr|\le(k-1)\sqrt p\]
where $\chi$ is the quadratic residue character (Legendre symbol) mod~$(p)$.

\medskip

\noindent{\bf 4.} For prime powers $q\equiv -1$  mod~$(4)$ the construction in \S\ref{Paley} yields the Paley tournament $T(q)$, a complete graph $K_q$ with directed edges, and the construction in \S\ref{InfPaley} yields, for each prime $p\equiv -1$ mod~$(4)$ and infinite set $E$ satisfying the conditions given there, an infinite Paley tournament $T(p^E)$. Again, there are uncountably many of these objects, but a slight adaptation of the preceding arguments shows that they are all isomorphic to the countable random tournament; a model of this can be obtained by applying the construction in \S\ref{Random} to primes $p, q\equiv -1$ mod~$(4)$, where quadratic reciprocity now gives $(\frac{p}{q})(\frac{q}{p})=-1$.

\medskip

\noindent{\bf 5.} Peter Cameron~\cite{Cam19} has suggested a more general construction using ultraproducts of finite fields, rather than direct limits, together with \L o\'s's Theorem, to approximate the random graph (see also~\cite[Example~1.3.6]{MS}, based on asymptotic classes and ultraproducts); this has the advantage of allowing finite fields of different characteristics to be used, thus yielding fields of characteristic $0$.


\section{Automorphism groups}\label{Auto}

It follows from a theorem of Carlitz~\cite{Car} that the automorphism group ${\rm Aut}\,P(q)$ of $P(q)$ is the subgroup $A\Delta L_1(q)$ of index $2$ in $A\Gamma L_1(q)$ consisting of the transformations
\[t\mapsto at^{\gamma}+b\quad (a, b\in{\mathbb F}_q,\, \chi(a)=1,\, \gamma\in{\rm Gal}\,{\mathbb F}_q)\]
of the vertex set ${\mathbb F}_q$, where ${\rm Gal}\,{\mathbb F}_q$ is the Galois group or automorphism group of ${\mathbb F}_q$, a cyclic group of order $\log_pq$ generated by the Frobenius automorphism $t\mapsto t^p$. The affine transformations (those elements with $\gamma=1$) and the translations (those with $\gamma=1$ and $a=1$) form normal subgroups $AHL_1(q)$ (`$H$' for `half') and $T_1(q)$ of $A\Delta L_1(q)$ with
\[A\Delta L_1(q)\ge AHL_1(q)>T_1(q)>1,\]
and the abelian quotients in this series show that $A\Delta L_1(q)$ is solvable, of derived length at most $3$.

One might hope that the automorphism group of $P(p^{rE})$ for $r=1$ or $2$ would have a similar structure. Clearly it contains the subgroup $A\Delta L_1(p^{rE})$ of index $2$ in $A\Gamma L_1(p^{rE})$ consisting of the transformations
\[t\mapsto at^{\gamma}+b\quad (a, b\in{\mathbb F}_{p^{rE}},\, \chi(a)=1,\, \gamma\in{\rm Gal}\,{\mathbb F}_{p^{rE}}).\]
Here ${\rm Gal}\,{\mathbb F}_{p^{rE}}$ is not the direct limit of the groups ${\rm Gal}\,{\mathbb F}_{p^{re}}$ for $e\in E$, but their {\sl inverse\/} limit: this can be identified with the (uncountable) subgroup of the cartesian product
$\prod_{e\in E}{\rm Gal}\,{\mathbb F}_{p^{re}}$
consisting of those elements whose coordinates $\gamma_{re}\in{\rm Gal}\,{\mathbb F}_{p^{re}}$ are consistent with the restriction mappings ${\rm Gal}\,{\mathbb F}_{p^{rf}}\to{\rm Gal}\,{\mathbb F}_{p^{re}}$ induced by inclusions ${\mathbb F}_{p^{re}}\subseteq{\mathbb F}_{p^{rf}}$ for $e$ dividing $f\in E$.

As in the finite case, this group $A\Delta L_1(p^{rE})$ is solvable, of derived length $3$. However, the facts that $P(p^{rE})\cong R$ and that ${\rm Aut}\,R$ acts transitively on isomorphism classes of finite induced subgraphs of $R$ (by the back-and-forth argument used in \S\ref{Random}) destroy any hope that this subgroup might be the whole of ${\rm Aut}\,P(p^{rE})$. Indeed, far from being solvable, ${\rm Aut}\,R$ has been shown by Truss~\cite{Tru} to be a simple group, and to contain a subgroup isomorphic to the symmetric group on a countably infinite set.

\section{Generalised Paley graphs}\label{GenPaley}

In 2009 Lim and Praeger~\cite{LP} introduced {\sl generalised Paley graphs\/} $P_d(q)$, where $q$ is a prime power $p^e$ and $d$ divides $q-1$ (for convenience, we have changed their notation). Again the vertex set is ${\mathbb F}_q$, but now vertices $x$ and $y$ are adjacent if and only if $x-y$ is contained in the unique subgroup $D$ of index $d$ in the multiplicative group ${\mathbb F}_q^*$, consisting of the non-zero $d$th powers. To give an undirected graph we assume that if $q$ is odd then the order $(q-1)/d$ of $D$ is even. For example, taking $d=2$ gives the Paley graphs $P(q)=P_2(q)$. 

The construction in \S\ref{InfPaley} carries through in the obvious way to give {\sl infinite generalised Paley graphs\/} $P_d(p^{rE})$ where $r$ is the multiplicative order of the prime $p$ mod~$(2d)$ (or mod~$(d)$ if $p=2$), except that we now need $E$ to consist of integers $e$ coprime to $d$. The proof of Theorem~\ref{PaleyRthm} also carries through, provided we take $\chi$ to be a multiplicative character of ${\mathbb F}_{q^e}$ of degree $d$ (equivalently with kernel $D$), and replace the factor $1+\chi(x-a)$ in equation~(\ref{abeqn}) with
\[1+\chi(x-a)+\chi(x-a)^2+\cdots+\chi(x-a)^{d-1}=\prod_{j=1}^{d-1}(\chi(x-a)-\omega^j)\]
where $\omega$ is a primitive $d$th root of $1$ in $\mathbb C$; again we can apply Weil's estimate, now in the more general form (\ref{Weil}) given in Remark~1, to show that $P_d(p^{rE})\cong R$.

The remarks in \S\ref{Auto} about automorphism groups also apply here, though it should be noted that, as shown in~\cite{LP}, there are examples where $d$ does not divide $p-1$ and ${\rm Aut}\,P_d(q)$ is significantly larger than the obvious analogue of $A\Delta L_1(q)$.

\section{Symmetry versus asymmetry}\label{Symm}

The main aim of Erd\H os and R\'enyi in~\cite{ER} was to consider, in the contexts of finite and countably infinite graphs, the balance between symmetric and asymmetric graphs, those with and without a non-identity automorphism. Most of the paper concerns finite graphs, and here they proved, in a very precise sense, that not only are most graphs asymmetric, but in fact they are on average a long way from being symmetric. For a finite graph $G=(V,E)$ they defined $A(G)$ to be the least number of edge-changes (insertions or deletions) required to convert $G$ into a symmetric graph on $V$. We may identify $G$ with its edge set $E$, regarded as an element of the power set ${\mathcal P}(V^{(2)})=({\mathbb F}_2)^{V^{(2)}}$ of the set $V^{(2)}$ of $2$-element subsets of $V$; the Hamming distance between two graphs $(V,E)$ and $(V,E')$, with respect to the basis consisting of the graphs with one edge, is $|E\oplus E'|$ where $\oplus$ denotes symmetric difference, so $A(G)$ is the distance from $G$ to the nearest symmetric graph on $V$.

For distinct vertices $u$ and $v$ in $G$ Erd\H os and R\'enyi defined $\Delta_{uv}$ to be the number of vertices $w\ne u, v$ adjacent to just one of $u$ and $v$. By making $\Delta_{uv}$ edge-changes one can give $u$ and $v$ the same neighbours, allowing an automorphism transposing them and fixing all other vertices, so
\[
A(G)\le \min_{u\ne v}\Delta_{uv}.
\]
By a simple counting argument they showed that if $G$ has order $n$ then
\begin{equation}\label{Deltabd}
\min_{u\ne v}\Delta_{uv}\le\lfloor\frac{n-1}{2}\rfloor,
\end{equation}
so that
\[
A(G)\le \lfloor\frac{n-1}{2}\rfloor.
\]
They then showed that `most' graphs $G$ of order $n$ have $A(G)$ close to $\lfloor(n-1)/2\rfloor$, so that they are very far from being symmetric. As an aside they defined a $\Delta$-{\sl graph\/} to be one achieving equality in (\ref{Deltabd}), and noted that the graphs $P(q)$ have this property: indeed, $\Delta_{uv}=(q-1)/2$ for all pairs $u\ne v$ in $P(q)$. Of course, these graphs are exceptional from this point of view, in that they satisfy $A(P(q))=0$.

By contrast, Erd\H os and R\'enyi showed in the last part of their paper that `most' countably infinite graphs are symmetric. Indeed, it follows from their construction of $R$ and the alternative one due to Rado~\cite{Rad} that most such graphs are isomorphic to $R$ and are therefore {\sl highly\/} symmetric: for example, ${\rm Aut}\,R$ acts transitively on finite induced subgraphs, and hence has rank~$3$ on the vertices. In fact, one can show that this group is uncountable, for instance by choosing a prime $p\equiv 1$ mod~$(4)$ and taking $E=\{q^n\mid n\ge 1\}$ in \S\ref{InfPaley} for some odd prime $q$, so that by our remarks in \S\ref{Auto} ${\rm Aut}\,R$ contains a copy of
\[{\rm Gal}\,P(p^E)=\lim_{\leftarrow}{\rm Gal}\,P(p^{q^n})\cong \lim_{\leftarrow}{\mathbb Z}/q^n{\mathbb Z}\cong {\mathbb Z}_q,\]
the uncountable group of $q$-adic integers.

\bigskip

\noindent{\bf Acknowledgment} The author is grateful to Peter Cameron and Dugald Macpherson for some very helpful comments.


\bigskip

\noindent School of Mathematics

\noindent University of Southampton

\noindent Southampton SO17 1BJ

\noindent UK

\medskip

\noindent G.A.Jones@maths.soton.ac.uk

\end{document}